\documentclass[12pt,a4paper, reqno]{amsart}
\usepackage{amssymb}
\usepackage{graphicx}
\usepackage{psfrag}

\newtheorem{theorem}{Theorem}[section]
\newtheorem{corollary}{Corollary}[section]

\newtheorem{definition}{Definition}[section]

\def \btheta {\mbox{\boldmath $\theta $}}

\def \Cov   {\text {\rm Cov}}

\def \Var   {\text {\rm Var}}
\def \diag  {\text {\rm diag}}

\def \Re    {\text {\rm Re}}
\def \ket    {\text {\it ket}}
\def \Ket    {\text {\rm ket}}

\def \diag  {\text {\rm diag}}
\def \lim   {\text {\rm lim}}

\def \tr    {\text {\rm tr}}
\def \tr    {\text {\rm Tr}}

\address{Indian Statistical Institute\\Delhi Centre, 7, S. J. S. Sansanwal Marg, New Delhi - 110 016, India\\e-mail: krp@isid.ac.in}
\begin{document}

%\maketitle

\begin{center}
{\large \bf On the philosophy of Cram\'er-Rao-Bhattacharya Inequalities in Quantum Statistics}
\vskip0.5in
{\it To my revered guru C. R. Rao for revealing the mystries of Chance }
\vskip0.5in
K. R. Parthasarathy\\
Indian Statistical Institute\\Delhi Centre, \\7, S. J. S. Sansanwal Marg,\\ New Delhi - 110 016, India\\e-mail: krp@isid.ac.in
\end{center}

\vskip0.3in
%\begin{abstract}
\noindent{\bf Summary \quad} To any parametric family of states of a finite level quantum system we associate a space of Fisher maps and introduce the natural notions of Cram\'er-Rao-Bhattacharya tensor and Fisher information form. This leads us to an abstract Cram\'er-Rao-Bhattacharya lower bound for the covariance matrix of any finite number of unbiased estimators of parameteric functions. A number of illustrative examples is included. Modulo technical assumptions of various kinds our methods can be applied to infinite level quantum systems as well as parametric families of classical probability distributions on Borel spaces.

\vskip0.1in
\noindent{\bf AMS Subject classification index:} 81C20, 94A15.

%\end{abstract}

\vskip0.9in
\section{Introduction}
The evolution of modern scientific thought is strewn with several examples expressing the following sentiment: in any effort to accomplish a task there can be a certain limit to the efficiency of its performance. In the present context we bring to attention three such famous examples which are based on the combination of a deep conceptual approach and simple mathematical arguments. Finally, we shall focus on one of them, namely, limits to the efficiency of estimating an unknown parameter involved in a family of states of a finite level quantum system.

Our first example is the celebrated uncertainty principle of Heisenberg \cite{wh} in quantum mechanics. For an interesting historical account of this great discovery in the philosophy of science we refer the reader to the essay by Jagdish Mehra in \cite{jm}. If $q$ and $p$ denote the position and momentum operators of a quantum mechanical particle executing motion on the real line $\mathbb{R}$ they obey the commutation relation $qp-pq = i \hslash$  where $\hslash = h/2 \pi,$ $h$ being the Planck's constant, and this implies the following inequality. If $\psi$ is the absolutely square integrable wave function describing the state of the system and $\Var (X | \psi)$  denotes the variance of the observable in the state $\psi$ then
\begin{equation}
\Var (q|\psi) \Var (p|\psi) \geqslant  \hslash^2/4. \label{eq1.1}
\end{equation}
In particular, if the variance of $p$ in the state $\psi$ is $\sigma^2$ then 
$$\Var (q|\psi) \geqslant \hslash^2/4 \sigma^2. $$
In other words, this sets a limit to the accuracy with which the position $q$ can be measured in the state $\psi.$ Such limits to accuracy hold for any `conjugate pair' of observables in quantum theory.

Our second example is the famous Cram\'er-Rao inequality \cite{hc}, \cite{crr} in the theory of estimation of statistical parameters. For an amusing and insightful account of the route by which this fundamental discovery was made and how it came to be recognized in the history of statistical science we refer to \cite{Crr}. Suppose $\{p (\omega, \theta)\}$ is a parametric family of probability density functions with respect to a $\sigma$-finite measure in a Borel space $(\Omega, \mathcal{F}),$ $\theta$ being a real parameter varying in an open interval $(a,b).$ Assume that the function
\begin{equation}
I (\theta) = \int_{\Omega} \left (\frac{\partial}{\partial \theta} \log p (\omega, \theta) \right )^2 p (\omega,\theta) \mu (d\omega)     \label{eq1.2}
\end{equation}
is well-defined for all $\theta$ in $(a,b).$ On the basis of a sample point $\omega$ obtained from experiment evaluate a function $T(\omega)$ as an estimate of the parameter $\theta.$ The function $T(\cdot)$ on $\Omega$ is called an {\it unbiased estimator} of $\theta$ if
$$\int_{\Omega} T (\omega) p (\omega, \theta) \mu(d\omega) = \theta \quad \forall \quad \theta \in (a,b)$$
and, in such a case, its variance, denoted by $V (T|\theta)$ is defined by
$$V (T|\theta) = \int (T(\omega)-\theta)^2 \,\,\, p (\omega, \theta) \mu (d\omega). $$
Indeed, $V(T | \theta)$ is a measure of the error involved in estimating $\theta$ by $T(\omega).$ The Cram\'er-Rao inequality in its simplest form says that
\begin{equation}
 V (T | \theta) \geqslant I (\theta)^{-1} \label{eq1.3}
\end{equation}
where $I(\theta)$ is given by \eqref{eq1.2} and called the `Fisher information' at $\theta.$ Thus \eqref{eq1.3} sets a limit to the accuracy of estimating the unknown parameter $\theta$ from experimental observation.

It is a remarkable fact that a special case of \eqref{eq1.3} implies the Heisenberg uncertainty principle \eqref{eq1.1} and much more in emphasizing the profundity of Fisher information. Indeed, let $\psi \in L^2 (\mathbb{R})$ be a  wave function so that $\|\psi\| = 1.$ By changing $\psi$ to a new wave function $e^{i \alpha x} \psi (x),$ $\alpha \in \mathbb{R},$ if necessary, we may assume, without loss of generality that the momentum operator $p$ satisfies the condition $\langle \psi | p | \psi \rangle = 0$ and $\langle \psi | q | \psi \rangle = m,$ a real scalar. By Born's interpretation $f = |\psi|^2$ is the probability density function of the position observable $q$ in the state $\psi.$ Introducing the parametric family $\{f (x-\theta), \theta \in \mathbb{R} \}$ of probability densities we see that its Fisher information $I (\theta)$ is given by
\begin{eqnarray}
I(\theta) &=&  \int_{\mathbb{R}} \left (\frac{f^{\prime} (x-\theta)}{f(x-\theta)} \right )^2 f (x-\theta) \,\,dx   \nonumber \\
&=&  4 \quad \int_{\mathbb{R}} \left ( \Re  \,\, \frac{\psi^{\prime}}{\psi} (x) \right )^2 \,\, | \psi (x) |^2 \,\,dx          \label{eq1.4}
\end{eqnarray}
and therefore independent of $\theta.$ By Cram\'er-Rao inequality we have
\begin{eqnarray}
\Var (q|\psi) &=& \int (x-m)^2 \,\, f (x) \,\,dx \nonumber \\
&\geqslant& \frac{1}{I(m)}. \label{eq1.5}
\end{eqnarray}
On the other hand
$$\Var (p|\psi) = \int_{\mathbb{R}} x^2 \,\,|(F \psi)(x)|^2 \,\,dx$$
where $F$ is the unitary Fourier transform in $L^2 (\mathbb{R}).$ Thus by \eqref{eq1.4} we have
\begin{eqnarray}
\Var (p|\psi) &=& \langle \psi | F^{\dagger} q^2 F | \psi \rangle \nonumber \\ 
&=& \|p \psi \|^2 \nonumber \\
&=& \hslash^2 \int \,\left |\psi^{\prime} (x) \right |^2 \,\,dx \nonumber \\
&=&\hslash^2 \int \,\left |\frac{\psi^{\prime}}{\psi} (x) \right |^2 \,\, \left |\psi (x) \right |^2\,\,dx \nonumber \\
&\geqslant& \hslash^2 \int \,\left |\left ( \Re \, \frac{\psi^{\prime}}{\psi}\right ) (x) \right |^2 \,\, \left |\psi (x) \right |^2\,\,dx \nonumber \\
&=& \frac{\hslash^2}{4} \,\,I(m) \label{eq1.6}
\end{eqnarray}
which together with \eqref{eq1.5} implies \eqref{eq1.1}. The more powerful inequality \eqref{eq1.6} and its natural generalization for covariance matrices in $L^2 (\mathbb{R}^n)$ are known together as Stam's uncertainty principle. For more information along these lines and a rich survey of information inequalities we refer to the paper \cite{dct} by A. Dembo, T. M. Cover and J. A. Thomas.

Our last illustrious example is of a different genre but again connected with the notion of information. It is Shannon's noisy coding theorem \cite{ceh} which sets a limit to the ability of communication through an information channel in the presence of noise. Again we present the simplest version of this strikingly beautiful result in order to highlight the philosophical aspect and refer to \cite{KrP} for more general versions. Consider an information channel whose input and output alphabets are same and equal to the binary alphabet $\{0,1\}$ which is also a field of two elements with the operations of addition and multiplication modulo 2. If an input letter $x$ from this alphabet is transmitted through the channel assume that the output letter is $x$ or $x+1$ with probability $1-p$ or $p$ so that the probability of error due to noise in transmission is $p.$ Such a channel is said to be {\it binary} and {\it symmetric.} Assume that the transmission of a sequence $x_1, x_2, \ldots, x_n$ of $n$ letters through this binary symmetric channel yields the output sequence $y_1, y_2, \ldots, y_n$ where $y_1 - x_1,$ $y_2 - x_2, \ldots, y_n - x_n$ are independently and identically distributed Bernoulli random variables, each assuming the values $0$ and $1$ with probability $q = 1-p$ and $p$ respectively. Such a channel is called a {\it memoryless} binary symmetric channel.

Denote the alphabet by $\mathbb{F}_2.$ By a {\it code} of {\it size} $m,$ {\it length} $n$ and {\it error probability} not exceeding $\varepsilon,$ where $0 < \varepsilon < 1,$ we mean $m$ pairs $(u_j, E_j), 1 \leqslant j \leqslant m,$ $u_j \in \mathbb{F}_2^n,$ $E_j \subset \mathbb{F}_2^n,$ $E_1, E_2, \ldots, E_m$ are pairwise disjoint, satisfying the inequalities
$$\mathbb{P} \,\left ( \mbox{output sequence} \,\, \in\,\, E_j | \,\mbox{input sequence} \,= u_j \right ) > 1 - \varepsilon \quad \forall \,\,j. $$
Denote by $N (n,p, \varepsilon)$ the maximum possible size for codes of length $n$ with error probability not exceeding $\varepsilon.$ Then
\begin{equation}
\underset{n \rightarrow \infty}{\lim}  \frac{1}{n} \,\log_2 \,N(n,p,\varepsilon) = 1+ p \,\log_2\, p + q \,\log_2 q \quad \forall \,\,0 < \varepsilon < 1, \,0 < p < 1/2.  \label{eq1.7}
\end{equation}
If we write $H(p) = - p \,\log_2 \,p - q \,\log_2 \,q$ and call it the Shannon entropy of the Bernoulli random variable with probability of success (error) $p$ then \eqref{eq1.7} has the interpretation that for large $n,$ among the $2^n$ possible input sequences of length $n$ roughly $2^{n(1-H(p))}$ sequences could be transmitted with error probability $< \,\varepsilon$ and not more. For this reason the expression on the right hand side of \eqref{eq1.7} is called the Shannon capacity of the binary symmetric channel with error probability $p.$ A corresponding generalization for memoryless and stationary quantum channels describing their `capacity' to transmit classical alphabetic messages exists. For a leisurely and self-contained exposition of such coding theorems see \cite{KrP}. The notion of entropy that arises in the brief discussion above can be introduced for a large class of density functions and this, in turn, leads to some remarkable connections with Fisher information and many powerful information theoretic inequalities. Once again we refer to the very rich survey article \cite{dct}.

All the three examples described above have been generalized in several ways, connections between them and relations with other branches of science and engineering have emerged and an enormous amount of literature has grown around them. The last example has given birth to the subject of quantum information theory and coding theorems for quantum channels \cite{man}, \cite{KrP}. The present essay is devoted to the second example but in the context of parametric families of states of finite level quantum systems. Starting from the books of Helstr\"om \cite{cwh}, Holevo \cite{ash}, and Hayashi \cite{mh} there is quite some literature on the Cram\'er Rao bounds for quantum systems. By confining ourselves to finite level systems we avoid the technical difficulties of dealing with unbounded operators and their varying domains but we gain conceptual and algebraic clarity.

In Section 2 we give a brief account of the quantum probability of finite level quantum systems in a complex finite dimensional Hilbert space including the notions of events, observables, states, generalized measurements and composite systems in the language of tensor products of Hilbert spaces. Heisenberg's uncertainty principle and an entropic uncertainty principle are briefly described. The notions of parametric families of states and unbiased estimators of parametric functions along with their variances and covariances are introduced.

Section 3 contains the key notions, namely,  Fisher maps, the Fisher information form and the Cram\'er-Rao-Bhattacharya (CRB) tensor with respect to a parametric family of states of a finite level quantum system. The  Cram\'er-Rao-Bhattacharya (CRB) bound is finally expressed in terms of the CRB tensor and the Fisher information form. Several illustrative examples are given.

In the last section we show how, by using a dilation theorem of Naimark, one can obtain a CRB bound for the covariance matrix of unbiased  estimators of parametric functions based on generalized measurements.

\section{Preliminaries in the quantum probability and statistics of finite level systems}
\setcounter{equation}{0}
\setcounter{theorem}{1}

A finite level quantum system is described by `states' in a finite dimensional complex Hilbert space. We choose and fix such a Hilbert space $\mathcal{H}$ with scalar product $\langle u|v \rangle $ which is linear in the variable $v$ and antilinear in $u.$ A typical example obtains when $\mathcal{H}$ is the $n$-dimensional complex vector space $\mathbb{C}^n$ of column vectors and its dual is the space of all row vectors. In this case the scalar product is expressed as
$$\langle u|v \rangle = \sum_{i=1}^n \bar{a}_i b_i $$
where
$$u = \left [ \begin{array}{c}  a_1 \\ a_2 \\ \vdots \\ a_n \end{array} \right ],  v = \left [ \begin{array}{c}  b_1 \\ b_2 \\ \vdots \\ b_n \end{array} \right ], \quad a_i, b_i \in \mathbb{C} \,\,\forall \,\, i.$$
Elements of $\mathcal{H}$ are called $\ket$ vectors, a typical element in $\mathcal{H}$ being denoted by $|v\rangle$  whereas any element in the dual of $\mathcal{H}$ is called a {\it bra} vector and a typical bra vector is denoted by $\langle u|.$ The linear functional $\langle u |$ evaluated at a $\Ket$ vector $|v\rangle$ is the scalar product $\langle u|v\rangle.$ If $A$ is an operator in $\mathcal{H}$ it is customary to write
$$\langle u |Av \rangle = \langle u|A|v \rangle. $$
The adjoint of $A$ is denoted by $A^{\dagger}$ so that
$$\langle u|A|v\rangle = \langle A^{\dagger} u | v \rangle = \langle u | Av \rangle.$$ 
In such a notation $|u\rangle \langle v|$ denotes the operator satisfying
$$\left ( |u \rangle \langle v|\right ) \,\, | w\rangle = \langle v|w\rangle | u \rangle \quad \forall \quad |u\rangle, |v\rangle, |w\rangle \,\,\mbox{in}\,\, \mathcal{H}.$$
The trace of an operator $A$ in $\mathcal{H}$ is denoted by $\tr \,A.$ In particular$\tr \, |u \rangle \langle v | = \langle v|u \rangle.$ Note that $|u \rangle \langle v |$ is a rank one operator when $|u \rangle \neq 0,$ $|v\rangle \neq 0,$ and
$$\left ( | u_1 \rangle \langle v_1 | \right ) \left ( | u_2\rangle \langle v_2 | \right ) \cdots \left ( |u_k \rangle \langle v_k | \right )= c |u_1 \rangle \langle v_k| $$
where $c = \langle v_1 | u_2 \rangle \langle v_2 | u_3 \rangle \cdots \langle v_{k-1} | u_k \rangle.$

Denote by $\mathcal{B}(\mathcal{H}),$ $\mathcal{P}(\mathcal{H}),$ $\mathcal{O}(\mathcal{H}),$ $\mathcal{S}(\mathcal{H})$ respectively the $\ast$- algebra of all operators on $\mathcal{H}$ with  its usual (strong) topology, the orthomodular lattice of all orthogonal projection operators on $\mathcal{H},$ the real linear space of all hermitian operators in $\mathcal{H}$ and the compact convex set of all nonnegative definite operators of unit trace. We have $\mathcal{P}(\mathcal{H}) \subset \mathcal{O}(\mathcal{H}) \subset \mathcal{B}(\mathcal{H})$ and $\mathcal{S}(\mathcal{H}) \subset \mathcal{O}(\mathcal{H}) \subset \mathcal{B}(\mathcal{H}).$ If $A, B \in \mathcal{O}(\mathcal{H})$ we say that $A \leqslant B$ if $B-A$ is nonnegative definite. Then $\mathcal{O}(\mathcal{H})$ is a partially ordered real linear space. A nonnegative definite hermitian operator is simply called a positive operator.

The zero and identity operators are denoted respectively by $O$ amd $I.$ Often, $I$ is denoted by $1.$ For any scalar $\lambda$ the operator $\lambda I$ is also denoted by $\lambda.$ Thus, for $A \in \mathcal{B}(\mathcal{H}),$ $\lambda \in \mathbb{C},$ $A - \lambda$ stands for the operator $A-\lambda I.$  For any $E \in \mathcal{P}(\mathcal{H}),$ $0 \leqslant E \leqslant 1$ and $(1-E) \in \mathcal{P}(\mathcal{H}).$ By a projection we shall always mean an orthogonal projection operator i.e., an element of $\mathcal{P}(\mathcal{H}).$ If $E_1, E_2 \in \mathcal{P}(\mathcal{H})$ and $E_1 \leqslant E_2$ then $(E_2-E_1) \in \mathcal{P}(\mathcal{H}).$ When a quantum system is described by $\mathcal{H}$ we say that the elements of $\mathcal{P}(\mathcal{H})$ are the {\it events} concerning the system, $0$ is the {\it null} event and $1$ is the {\it certain} event. If $E_1, E_2 \in \mathcal{P}(\mathcal{H})$ and $E_1 \leqslant E_2$ we say that the event $E_1$ {\it implies} the event $E_2.$ If $E \in \mathcal{P}(\mathcal{H})$ then $1-E$ is the event `not \,$E$'. If $E_1, E_2 \in \mathcal{P}(\mathcal{H})$ their maximum $E_1 \vee E_2$ and minimum $E_1 \wedge E_2$ are respectively interpreted as `$E_1 \,\mbox{or}\, E_2$' and `$E_1 \,\mbox{and}\,E_2$'. If $E_1 E_2 = 0$ then $E_1 \vee E_2 = E_1 + E_2.$ If $E_1$ and $E_2$ commute then $E_1 \wedge E_2 = E_1 E_2.$ The first basic difference between classical probability and quantum probability theory arises from the fact that for three events $E_i$ in $\mathcal{P}(\mathcal{H}), i=1,2,3$ one may not have $E_1 \wedge (E_2 \vee E_3) = (E_1 \wedge E_2) \vee (E_1 \wedge E_3).$ Whenever the $E_i$'s commute with each other the operations $\wedge$ and $\vee$ distribute with each other.

Any hermitian opearator $X$ in $\mathcal{H}$ is called a real-valued or simply an {\it observable} about the system described by $\mathcal{H}.$ Thus $\mathcal{O}(\mathcal{H})$ is the real linear space of all real-valued observables. If $X, Y \in \mathcal{O}(\mathcal{H})$ and $XY=YX$ then $XY$ is also an element of $\mathcal{O}(\mathcal{H}).$ If $X \in \mathcal{O}(\mathcal{H})$ and $\sigma(X)$ denotes the set of all its eigenvalues then, by the spectral theorem, $X$ admits a unique spectral resolution or representation
\begin{equation}
X = \sum_{\lambda \in \sigma (X)} \lambda \, E_{\lambda}  \label{eq2.1}
\end{equation}
where $\sigma (X) \subset \mathbb{R}$ is a finite set of cardinality not exceeding the dimension of $\mathcal{H},$ $0 \neq E_{\lambda} \in \mathcal{P}(\mathcal{H})$ $\forall$ $\lambda \in \sigma (X)$ and
\begin{eqnarray}
\sum_{\lambda \in \sigma (X)} \,E_{\lambda} &=&  I,  \label{eq2.2} \\
E_{\lambda} E_{\lambda^{\prime}} &=& \delta_{\lambda \lambda^{\prime}} E_{\lambda} \quad \forall \quad \lambda, \lambda^{\prime} \in \sigma (X).\label{eq2.3}
\end{eqnarray}
This, at once, suggests the interpretation that the eigenprojection $E_{\lambda}$ associated with the eigenvalue $\lambda$ in \eqref{eq2.1} is the event that the observable $X$ takes the value $\lambda$ and $\sigma(X)$ is the set of all values that $X$ can take. If $\varphi : \sigma (X) \rightarrow \mathbb{R}$ or $\mathbb{C}$ is a real or complex-valued function then
\begin{equation}
\varphi (X) = \sum_{\lambda \in \sigma (X)} \,\,\varphi(\lambda) \,E_{\lambda}   \label{eq2.4}
\end{equation}
is the real or complex-valued observable which is the function $\varphi$ of $X.$

Any element $\rho \in \mathcal{S}(\mathcal{H})$ is called a {\it state} of the quantum system described by $X.$ Such a state $\rho$ is also called a {\it density operator.} Clearly, $\rho$ itself becomes an observable. If $E \in \mathcal{P}(\mathcal{H})$  is an event and $\rho$ is a state then $\tr \, \rho E$ is a quantity in the unit interval $[0,1]$ called the {\it probability} of the event $E$ in the state $\rho.$ If $E_1, E_2$ are two events satisfying the relation $E_1 E_2 = 0$ then $E_1 + E_2$ is also an event and $\tr \, \rho (E_1 + E_2) = \tr \,\rho E_1 + \tr \,\rho E_2.$ However, for two events, $E_1, E_2$ it is not necessary that $\tr \,\rho (E_1 \vee E_2) \leqslant \tr \,\rho E_1 + \tr \, \rho E_2.$ In short, subadditivity property for probability need not hold good. But this property is retained whenever $E_1$ and $E_2$ commute with each other.

If $\rho$ is a state and $X$ is an element of $\mathcal{O}(\mathcal{H})$ with spectral resolution \eqref{eq2.1} then $\tr \,\rho E_{\lambda}$ is the probability that $X$ takes the value $\lambda$ in the state $\rho$ whenever $\lambda \in \sigma (X).$ Thus the {\it expectation} of $X$ in the state $\rho$ is given by
$$\sum_{\lambda \in \sigma (X)} \,\, \lambda \,\, \tr \,\rho E_{\lambda} = \tr \, \rho \sum_{\lambda \in \sigma (X)} \,\, \lambda E_{\lambda} = \tr \,\rho X.$$
More generally, the expectation of $\varphi (X)$ defined by \eqref{eq2.4} is given by $\tr \,\rho \varphi (X).$ In particular, the {\it variance} of $X$ in the state $\rho,$ denoted by $\Var (X|\rho)$ is given by
\begin{equation*}
 \begin{split}
\Var (X|\rho) &= \tr \, \rho X^2 - (\tr \,\rho X)^2 \\
&= \tr \, \rho (X - m)^2
 \end{split}
\end{equation*}
where $m = \tr \, \rho X$ is the expectation or {\it mean} of $X$ in the state $\rho.$ This shows, in particular, that $\Var (X|\rho)$ vanishes if and only if the restriction of $X$ to the range of $\rho$ is a scalar multiple of the identity.

The spectral theorem implies that the extreme points of the convex set $\mathcal{S}(\mathcal{H})$ are one dimensional projections of the form $|\psi\rangle \langle \psi |$ where $| \psi \rangle$ is a unit vector in $\mathcal{H}.$ Here, the projection remains unaltered if $| \psi \rangle$ is replaced by $c | \psi \rangle$ where $c$ is a scalar of modulus unity. Extreme points of $\mathcal{S}(\mathcal{H})$ are called {\it pure states} and a pure state is a one dimensional projection which, in turn, is determined by a unit vector in $\mathcal{H}$ modulo a scalar of modulus unity. By abuse of language any determining unit vector itself is called a pure state. Thus whenever we say that a unit vector $| \psi \rangle$ is a pure state we mean the density operator $| \psi \rangle \langle \psi |.$ By spectral theorem any state $\rho$ can be expressed as $\sum\limits_j p_j | \psi_j \rangle \langle \psi_j |$ where $p_1, p_2, \ldots$ is a finite probability distribution and $\left \{| \psi_j \rangle, j = 1,2, \ldots  \right \}$ is an orthonormal set of vectors in $\mathcal{H}.$ If $\{|\psi_j \rangle\}$ is any set of unit vectors and $p_j, j =1,2, \ldots$ is a probability distribution then $\sum\limits_j p_j | \psi_j \rangle \langle \psi_j |$ is a state. If $| \psi \rangle$ is a pure state and $X$ is a real-valued observable then its variance $\Var (X| |\psi \rangle)$ in the pure state $|\psi \rangle$ is zero if and only if $|\psi \rangle$ is an eigenvector for $X.$ Thus, even in a pure state $| \psi \rangle,$ an observable need not have a degenerate distribution. This is a significant departure from classical probability.

Hereafter, unless otherwise explicitly mentioned, we shall mean by an observable a real-valued observable. Let $X,Y$ be two observables, $\rho$ a state and let $m = \tr \,\rho X,$ $m^{\prime} = \tr \, \rho Y$ their respective means. Put $\widetilde{X} = X - m,$ $\widetilde{Y} = Y - m^{\prime}$ and consider the nonnegative function
$$f(z) = \tr \, \rho(\widetilde{X} + z \widetilde{Y})^{\dagger} (\widetilde{X} + z \widetilde{Y}), \quad z \in \mathbb{C}. $$
Then the inequality
$$\inf_{z \in \mathbb{C}} f (z) \geqslant 0$$
implies (see \cite{cwg}, \cite{KRP})
\begin{equation}
\Var (X|\rho) \Var (Y|\rho) \geqslant \left \{ \tr \,\rho \,\frac{1}{2i} (\widetilde{X} \widetilde{Y} - \widetilde{Y} \widetilde{X}) \right \}^2 + \left \{\tr \, \rho \, \frac{1}{2}  (\widetilde{X} \widetilde{Y} - \widetilde{Y} \widetilde{X})\right \}^2   \label{eq2.5}
\end{equation}
and thus puts a lower bound on the product of the variances of $X$ and $Y$ in a state $\rho.$ The quantum probability of finite level systems we have described here has a natural generalization when $\mathcal{H}$ is an infinite dimensional Hilbert space. When $\mathcal{H}=L^2 (\mathbb{R}),$ $X = q,$  $Y= p$ are the well-known position and momentum operators satisfying the Heisenberg commutation relations $qp-pq = i \hslash$ the inequality \eqref{eq2.5} yields the special form
$$\Var (q \big | | \psi \rangle) \Var (p \big | | \psi \rangle) \geqslant \frac{\hslash^2}{4}  \quad \forall \quad | \psi \rangle \in \mathcal{D}$$
where $\mathcal{D}$ is a dense domain in $\mathcal{H}$ where unbounded operators like $qp, pq$ etc. are well-defined. Thus \eqref{eq2.5} is at the heart of the Heisenberg's principle of uncertainty.

Now we introduce a notion which is  more general than that of an observable. Indeed, it plays an important role in the quantum version of Shannon's coding theorems of classical information theory.

\vskip0.1in
\begin{definition} \label{def2.1}
A  generalized measurement $\mathcal{L}$ of a finite level quantum system with Hilbert space $\mathcal{H}$ is a pair $(S, L)$ where $S$ is a finite set and $L : S \rightarrow \mathcal{B}(\mathcal{H})$ is a map satisfying the condition:
\begin{equation}
\sum_{s \in S} L (s)^{\dagger} L(s) = I.\label{eq2.6}
\end{equation}
\end{definition}

Such a generalized measurement $\mathcal{L} = (S,L)$ has the following interpretation. If the system is in the state $\rho$ and the measurement $\mathcal{L}$ is performed then the `value' $s \in S$ is obtained with probability $\tr \, L(s) \rho L(s)^{\dagger}$ and the system `collapses' to a new state
\begin{equation}
  \frac{L(s) \rho L(s)^{\dagger}}{\tr \, L(s) \rho L(s)^{\dagger}}. \label{eq2.7}
\end{equation}
If, for example, the system is initially in the state $\rho,$ a generalized measurement $\mathcal{L}_1 = (S_1, L_1)$ is performed and followed by another generalized measurement $\mathcal{L}_2 = (S_2, L_2)$ then the probability of obtaining the value $s_1 \in S_1$ is $\tr \, L_1 (s_1) \rho L_1 (s_1)^{\dagger}$ and the conditional probability of getting the value $s_2 \in S_2$ from $\mathcal{L}_2$ given the value $s_1$ from $\mathcal{L}_1$ is 
$$\tr \, L_2 (s_2) \left \{\frac{L_1 (s_1) \rho L_1 (s_1)^{\dagger}}{\tr \, L_1 (s_1) \rho L_1 (s_1)^{\dagger}} \right \} L_2 (s_2)^{\dagger}. $$
Thus the probability of obtaining the value $(s_1, s_2)$ from $\mathcal{L}_1$ followed by $\mathcal{L}_2$ is equal to
$$p(s_1, s_2) = \tr \, L_2 (s_2) L_1 (s_1) \rho \, L_1 (s_1)^{\dagger} L_2 (s_2)^{\dagger} $$
and the final collapsed state is
$$\frac{L_2 (s_2) L_1(s_1) \rho L_1 (s_1)^{\dagger} L_2 (s_2)^{\dagger}}{p (s_1, s_2)}.$$
More generally, if the measurements $\mathcal{L}_i = (S_i, L_i),$ $i =1,2,\ldots,m$ are performed in succession on a quantum system with initial state $\rho$ then the probability \linebreak $p(s_1, s_2, \ldots, s_m)$ of getting the sequence $s_1, s_2, \ldots, s_m$ of values $s_j \in S_j \,\forall \,\,j$ is given by
$$p (s_1, s_2, \ldots, s_m) = \tr \, L_m (s_m) L_{m-1} (s_{m-1}) \cdots L_1 (s_1) \rho L_1 (s_1)^{\dagger} L_2 (s_2)^{\dagger} \cdots L_m (s_m)^{\dagger} $$
and the final collapsed state is
$$\frac{1}{p(s_1, s_2, \ldots, s_m)} L_m (s_m) L_{m-1} (s_{m-1}) \ldots L_1 (s_1) \rho L_1 (s_1)^{\dagger} L_2 (s_2)^{\dagger} \ldots L_m (s_m)^{\dagger}.$$
This at once suggests the product rule for measurements $\mathcal{L}_i = (S_i, L_i)$ $i=1,2$ as $\mathcal{L} = (S_1 \times S_2, \widetilde{L})$ where
$$\widetilde{L} (s_1, s_2) = L_2 (s_2) L_1 (s_1), \quad s_1 \in S_1, s_2 \in S_2. $$
The measurement $\mathcal{L}$ stands for the measurement $\mathcal{L}_1$ followed by the measurement $\mathcal{L}_2.$

If $\mathcal{L} = (S,L)$ is a measurement with $S \subset \mathbb{R}$ or $\mathbb{C}$ then its {\it expectation} in the state $\rho$ is given by
$$\sum_{s \in S} s \,\tr\, L(s) \rho L(s)^{\dagger} = \sum_{s \in S} \,s\,\tr\,\rho L(s)^{\dagger} L(s). $$
If $S \subset \mathbb{R}$ its variance $\Var (\mathcal{L}|\rho)$ in the state $\rho$ is given by
$$\sum_{s \in S} s^2 \,\tr \, \rho \, L(s)^{\dagger} L(s) - \left ( \sum_{s \in s}\,s \,\tr\, \rho L(s)^{\dagger} L(s) \right )^2. $$

When $L(s)$ is a projection for every $s \in S$ then $(S,L)$ is called a {\it projective} or {\it von Neumann} {\it measurement}. If, in addition, $S \subset \mathbb{R}$ then the hermitian operator $\sum\limits_{s \in S} \,s \,L(s)$  is an observable and our notion of generalized measurement reduces to measuring an observable. It may be of some interest to formulate and obtain an uncertainty principle for a pair of two real-valued measurements.

For a measurement with values in an abstract set $S$ it is natural to replace the notion of variance by its entropy in a state $\rho.$ Thus we consider the quantity
\begin{equation}
H(\mathcal{L}|\rho) = - \sum_{s \in S} \,p(s)\, \log_2 \,p(s) \label{eq2.8}
\end{equation}
where
$$p(s) = \,\tr \, \rho \,L(s)^{\dagger}\, L(s) $$
and call it the {\it entropy} of the measurement $\mathcal{L} = (S, L)$ in the state $\rho.$ With this definition one has the following entropic uncertainty principle.
\vskip0.1in
\begin{theorem}[\cite{mkkrp}, \cite{hmjbmu}] \label{thm2.2}
Let $\mathcal{L} = (S,L),$ $\mathcal{M} = (T,M)$ be two generalized measurements of a finite level quantum system in a Hilbert space $\mathcal{H}.$ Let $L(s)^{\dagger} L(s) = X(s),$ $M(t)^{\dagger} M(t) = Y(t),$ $s \in S,$ $t \in T.$ Then for any state $\rho$ the following holds:
\begin{equation}
H (\mathcal{L}|\rho) + H(\mathcal{M}|\rho) \geqslant - 2 \log_2 \,\underset{s,t}{\max} \, \big| \big |X(s)^{1/2} Y(t)^{1/2} \big | \big |. \label{eq2.9}
\end{equation}
\end{theorem}
\vskip0.1in
\noindent{\bf Remark \quad} It is important to note that the right hand side in the inequality \eqref{eq2.9} is independent of $\rho.$

\vskip0.2in
If $X_i,$ $1 \leqslant i \leqslant k$ are $k$ observables, $\rho$ is a state in $\mathcal{H}$ and $\tr \,\rho X_i = m_i$ define the scalar
\begin{equation}
\nu_{ij} = \frac{1}{2} \,\tr \, \rho \left \{ (X_i - m_i)(X_j - m_j)+(X_j-m_j)(X_i - m_i)\right \}.   \label{eq2.10}
\end{equation}
Then the real symmetric matrix $((\nu_{ij}))$ of order $k$ is called the {\it covariance matrix} of the observables $X_1, X_2, \ldots, X_k$ in the state $\rho$ and denoted by $\Cov \left (X_1, X_2, \ldots, X_k \big | \rho \right ).$ It is a positive semidefinite matrix and it is  important to note the symmetrization in $i,j$ in the right hand side of \eqref{eq2.10}. Without such a symmetrization $\nu_{ij}$ could be a complex scalar.

Till now we talked about a single quantum system. Suppose we have a composite quantum system made out of several simple systems $A_1, A_2, \ldots, A_k$ with respective Hilbert spaces $\mathcal{H}_{A_{1}},\mathcal{H}_{A_{2}},  \ldots, \mathcal{H}_{A_{k}}.$ Then the Hilbert space of the joint system $A_1 A_2 \ldots A_k$ is the tensor product
$$\mathcal{H}_{A_{1} \ldots A_{k}}  = \mathcal{H}_{A_{1}} \otimes  \mathcal{H}_{A_{2}} \otimes \cdots \otimes   \mathcal{H}_{A_{k}}.$$
This is the quantum probabilistic analogue of cartesian product of sample spaces in classical probability. It is clear that
$$\dim \, \mathcal{H}_{A_{1} \ldots A_{k}} = \prod_{i=1}^k \dim \, \mathcal{H}_{A_{i}},$$
$\dim$ indicating dimension. If $\rho_i$ is a state in $\mathcal{H}_{A_{i}} \,\forall \, i$ then $\rho_1 \otimes \cdots \otimes \rho_k$ is a state of the composite system $A_1 A_2 \ldots A_k$ called the {\it product} state. If $\rho$ is a state in $\mathcal{H}_{A_{1}\ldots A_{k}}$ and  we take its relative trace over $\mathcal{H}_{A_{i_{1}}}, \mathcal{H}_{A_{i_{2}}}, \ldots, \mathcal{H}_{A_{i_{\ell}}}$ then we get the marginal state of the system $A_{r_{1}}, A_{r_{2}} \ldots, A_{r_{m}}$ where $\{1,2, \ldots, k\}$ is the disjoint union $\{i_1, i_2, \ldots, i_{\ell} \} \cup \{r_1, r_2, \ldots, r_m \}$ with $\ell + m = k.$ In this context of composite quantum systems there arises a new distinguishing feature of the subject with a remarkable role in physics as well as information theory. It is the existence of a very rich class of states in $\mathcal{H}_{A_{1} A_{2} \ldots A_{k}}$ which do not belong to the convex hull of all product states. Such states are called {\it entangled states} and they constitute a rich resource in quantum communication \cite{man}.

Till now we restricted ourselves to quantum probability. Now we describe a few basic concepts in quantum statistics dealing with a parametric family of quantum states of a finite level system. Let $\Gamma$ be a {\it parameter space} and let $\{\rho (\theta), \theta \in \Gamma \}$ be a {\it parametric family} of states in a Hilbert space $\mathcal{H}.$ Suppose $X$ is an observable, i.e., an element of $\mathcal{O}(\mathcal{H})$ and
\begin{equation}
\tr \,\,\rho (\theta) X = f (\theta), \,\, \theta \in \Gamma, \label{eq2.11}
\end{equation}
where $f$ is a real-valued function on $\Gamma.$ then we say that the observable $X$ is an {\it unbiased estimator} of the {\it parametric function} $f$ on $\Gamma.$

When the parametric family $\{ \rho (\theta), \theta \in \Gamma \}$ is fixed we write
\begin{equation}
\Var (X | \theta) = \Var (X|\rho(\theta)) \label{eq2.12}
\end{equation}
If $X_1, X_2, \ldots, X_m$ are $m$ observables we write
\begin{equation}
\Cov (X_1, \ldots, X_m | \theta) = \Cov (X_1, \ldots, X_m | \rho (\theta)). \label{eq2.13}
\end{equation}
A real-valued function $f$ on $\Gamma$ is said to be {\it estimable} with respect to $\{\rho(\theta), \theta \in \Gamma \}$ if there exists an observable $X \in \mathcal{O}(\mathcal{H})$ such that
$$\tr \,\rho (\theta) X = f (\theta) \,\,\forall \,\, \theta \in \Gamma.$$
Denote by $\mathcal{E}(\Gamma)$ the real linear space of all such estimable functions. An observable $X$ is said to be {\it balanced} with respect to the family $\{\rho (\theta), \theta \in \Gamma\}$ if $\tr \, \rho (\theta) X = 0 \,\,\forall \,\, \theta \in \Gamma.$  Denote by $\mathcal{N}$ the real linear space of all such balanced observables. For any $f \in \mathcal{E} (\Gamma),$ an unbiased estimator $X$ of $f$ write
$$\nu_f (\theta) = \inf \left \{\Var (X+Z|\theta), \,\,Z \in \mathcal{N} \right \}. $$
It is natural to look for good lower bounds for the function $\nu_f (\theta).$ We shall examine this problem in the next section and study some examples. If $f_j, \,1 \leqslant j \leqslant m$ are estimable parametric functions we shall also look for matrix lower bounds for the positive semidefinite matrices $\Cov (X_1, \ldots, X_m | \theta)$ as each $X_i$ varies over all unbiased estimators of $f_i$ for each $i =1,2, \ldots, m.$

For a more detailed introduction to quantum probability theory we refer to \cite{Krp}, \cite{KRP}. For an initiation to estimation theory and testing hypotheses in quantum statistics we refer to \cite{mh}, \cite{cwh}, \cite{ash}, References \cite{mh}, \cite{ash}, \cite{man}, \cite{KrP} contain applications of the theory of generalized measurements.

\section{The Fisher information form and the Cram\'er-Rao-Bhattacharya tensor} 
\setcounter{equation}{0}
\setcounter{definition}{0}

We consider a fixed parametric family $\{ \rho (\theta), \theta \in \Gamma \}$ of states of a finite level quantum system in a Hilbert space $\mathcal{H}$ with parameter space $\Gamma.$ As mentioned in the preceding section denote by $\mathcal{E}(\Gamma)$ and $\mathcal{N}$ respectively the real linear spaces of estimable functions and balanced observables. Recall that for any two unbiased estimators $X$ and $Y$ of an element $f \in \mathcal{E}(\Gamma),$ the observable $X-Y$ is an element of $\mathcal{N}.$

\vskip0.1in
\begin{definition}\label{def3.1}
A map $F : \Gamma \rightarrow \mathcal{B}(\mathcal{H})$ is called a {\it Fisher map} for the family $\{\rho (\theta), \theta \in \Gamma\}$ of states in $\mathcal{H}$ if the following two conditions hold:
\begin{itemize}
 \item[(i)] $\tr \, \rho (\theta) F(\theta) = 0 \quad \forall \,\,\theta \in \Gamma,$
\item[(ii)] $\tr \, \rho (\theta) \left \{F(\theta)^{\dagger} X + X F(\theta) \right \} = 0 \quad \forall \,\,\theta \in \Gamma, X \in \mathcal{N}.$
\end{itemize}
\end{definition}

Denote by $\mathcal{F}$ the real linear space of all Fisher maps with respect to $\{ \rho(\theta), \theta \in \Gamma \}$ and by $\mathcal{A}(\Gamma)$ the algebra of all real-valued functions on $\Gamma.$ If $a \in \mathcal{A}(\Gamma)$ and $F \in  \mathcal{F}$ then $aF$ defined by  $(aF) (\theta) = a(\theta) F(\theta)$ is also in $\mathcal{F}.$ In other words $\mathcal{F}$ is an $\mathcal{A}(\Gamma)$-module. For any two Fisher maps $F,G$ in $\mathcal{F}$ define
\begin{equation}
\begin{split}
\mathcal{I} (F,G)(\theta) & = \tr \, \rho (\theta) \frac{1}{2} \left (F(\theta)^{\dagger}  G(\theta) + G(\theta)^{\dagger} F(\theta) \right ) \\
&= \Re \,\tr \, \rho(\theta) F (\theta)^{\dagger} G(\theta). \label{eq3.1}
\end{split}
\end{equation}
Then $\mathcal{I}$ is called the {\it Fisher information form} associated with $\{\rho (\theta), \theta \in \Gamma \}.$ It may be noted that, for all $F, F_1, F_2, G \in \mathcal{F}$ and $a \in \mathcal{A} (\Gamma),$ the following hold:
\begin{equation*}
\begin{split}
\mathcal{I} (F,G) &= \mathcal{I} (G,F), \\
\mathcal{I} (aF,G) & = a \,\mathcal{I} (F,G), \\
\mathcal{I} (F_1 + F_2, G) &= \mathcal{I}(F_1,G) + \mathcal{I} (F_2,G),\\
\mathcal{I}(F,F) & \geqslant 0.
\end{split}
\end{equation*}
In particular, for any $F_i,$ $1 \leqslant i \leqslant n$ in $\mathcal{F}$ the matrix
\begin{equation}
\mathcal{I}_n (F_1, F_2, \ldots, F_n | \theta) = \left (\left (\mathcal{I} (F_i, F_j) (\theta) \right ) \right ), \theta \in \Gamma, i,j \in \{ 1,2, \ldots, n\} \label{eq3.2}
\end{equation}
is positive semidefinite. It is called the {\it information matrix} at $\theta$ corresponding to the elements $F_i,$ $ 1 \leqslant i \leqslant n$ in $\mathcal{F}.$

If $f \in \mathcal{E} (\Gamma),$ $ F \in \mathcal{F}$ define
\begin{equation}
\lambda (f,F) (\theta) = \tr \, \rho (\theta) \frac{1}{2} \left ( F(\theta)^{\dagger} X + X F (\theta) \right ), \quad \theta \in \Gamma  \label{eq3.3}
\end{equation}
where $X$ is any unbiased estimator of $f.$  Note that, in view of property (ii) in Definition \ref{def3.1} the right hand side of \eqref{eq3.3} is independent of the choice of the unbiased estimator of $f.$ Clearly, $\lambda (f, F)$ is real linear in the variable $f$ when $F$ is fixed and $\mathcal{A} (\Gamma)$-linear in the variable $F$ when $f$ is fixed. Thus $\lambda (\cdot , \cdot)$ can be viewed as an element of $\mathcal{E}(\Gamma) \otimes \mathcal{F}.$ We call $\lambda (\cdot , \cdot)$ the {\it Cram\'er-Rao-Bhattacharya tensor} or simply the  {\it CRB-tensor} associated with $\{\rho(\theta), \theta \in \Gamma\}.$

Let $f_i \in  \mathcal{E} (\Gamma),$ $X_i$ an unbiased estimator of $f_i$ for each $1 \leqslant i \leqslant m$ and let $F_j,$ $1 \leqslant j \leqslant n$ be Fisher maps with respect to $\{ \rho (\theta), \theta \in \Gamma \}.$ Define the $m \times m$ matrix
\begin{eqnarray}
\Lambda_{mn} (\theta) & =& ((\lambda_{ij} (\theta))), \quad 1 \leqslant i \leqslant m, \quad 1 \leqslant j \leqslant n, \quad \theta \in \Gamma \label{eq3.4}\\
\lambda_{ij} (\theta) &=& \lambda (f_i, F_j) (\theta) \quad \theta \in \Gamma, \label{eq3.5}
\end{eqnarray}
$\lambda$ being the CRB-tensor. We now introduce the family of positive semidefinite sesquilinear forms indexed by $\theta \in \Gamma$ in the vector space $\mathcal{B}(\mathcal{H})$ by
\begin{equation}
B_{\theta}  (X,Y) = \tr \, X^{\dagger} \,\rho (\theta) Y, \quad X,Y \in \mathcal{B}(\mathcal{H}). \label{eq3.6}
\end{equation}
By property (i) in Definition \ref{def3.1}, equations \eqref{eq3.3} and \eqref{eq3.5} we have
\begin{eqnarray*}
\lambda_{ij} (\theta) &=& \tr \, \rho (\theta) \frac{1}{2} \left \{ F_j(\theta)^{\dagger} (X_i - f_i (\theta)) + (X_i - f_i (\theta)) F_j (\theta)\right  \} \\
&=& \Re \, B_{\theta} \left ( X_i - f_i (\theta), F_j (\theta)^{\dagger} \right ).
\end{eqnarray*}
Multiplying both sides by the real scalars $a_i b_j$ and adding over $1 \leqslant i \leqslant m,$ $1 \leqslant i \leqslant n,$ we obtain
\begin{equation}
\mathbf{a}^{\prime} \Lambda_{mn} (\theta) \mathbf{b} = \Re \,\, B_{\theta} \left (\sum_{i=1}^m a_i (X_i - f_i (\theta)), \sum_{j=1}^n b_j F_j (\theta)^{\dagger}
 \right )   \label{eq3.7}
\end{equation}
where $\Lambda_{mn}$ and $B_{\theta}$ are given by \eqref{eq3.4}, \eqref{eq3.5} and \eqref{eq3.6} and $\mathbf{a}, \mathbf{b}$ are respectively column vectors of length $m,n$ with prime $^{\prime}$ indicating transpose. Now an application of the Cauchy-Schwarz inequality to the right hand side of \eqref{eq3.7} implies
\begin{eqnarray*}
\left ( \mathbf{a}^{\prime} \Lambda_{mn} (\btheta) \,\mathbf{b} \right )^2   &\leqslant& B_{\theta} \left (\sum_{i=1}^n a_i (X_i - f_i (\theta)), \sum_{i=1}^m a_i (X_i - f_i (\theta))  \right ) \\
&& \times B_{\theta} \left (\sum_{j=1}^n b_j F_j (\theta)^{\dagger}, \sum_{j=1}^n b_j F_j (\theta)^{\dagger}   \right ) \\
&=& \left \{\mathbf{a}^{\prime} \Cov \left ( X_1, X_2, \ldots, X_m | \theta \right ) \mathbf{a} \right \} \left \{\mathbf{b}^{\prime} \mathcal{I}_n \left (F_1, F_2, \ldots, F_n | \theta \right ) \mathbf{b} \right \}.
\end{eqnarray*}
Dividing both sides of this inequality by $\mathbf{b}^{\prime} \mathcal{I}_n \left (F_1, F_2, \ldots, F_n | \theta \right ) \mathbf{b} ,$ fixing $\mathbf{a}$ and maximizing the left hand side over all $\mathbf{b}$ satisfying $\mathcal{I}_n \left (F_1, F_2, \ldots, F_n | \theta \right ) \mathbf{b} \neq \mathbf{0}$ we obtain the matrix inequality:
$$\Lambda_{mn} (\theta) \mathcal{I}_n^{-} (F_1, F_2, \ldots, F_n | \theta) \Lambda_{mn} (\theta)^{\prime} \leqslant \Cov (X_1, X_2, \ldots X_m | \theta), $$
$\mathcal{I}_n^{-}$ denoting the generalized inverse of $\mathcal{I}_n (F_1, F_2, \ldots, F_n | \theta).$ In other words we have proved the following theorem
\vskip0.2in
\begin{theorem}[Quantum Cram\'er-Rao-Bhattacharya (CRB) inequality] \label{thm3.1}
Let \linebreak $\{\rho (\theta), \theta \in \Gamma\}$ be a parametric family of states of a finite level quantum system in a Hilbert space $\mathcal{H},$ $f_i,$ $1 \leqslant i \leqslant m$ estimable functions on $\Gamma,$ $X_i$ an unbiased estimator of $f_i$ for each $i$ and let $F_j,$ $ 1 \leqslant j \leqslant n$ be Fisher maps with respect to $\{\rho (\theta), \theta \in \Gamma \}.$ Then the following matrix inequality holds:
$$\Cov \left ( X_1, X_2, \ldots, X_m| \theta \right ) \geqslant \Lambda_{mn} (\theta) \mathcal{I}_n^{-} \left ( F_1, F_2, \ldots, F_n | \theta \right ) \Lambda_{mn} (\theta)^{\prime} \,\,\forall \,\, \theta \in \Gamma$$
where $\Lambda_{mn} (\theta)$ is the $m \times n$ matrix defined by \eqref{eq3.3}-\eqref{eq3.5} and $\mathcal{I}_n^{-} (F_1, F_2, \ldots, F_n | \theta)$ is the generalized inverse of the Fisher information matrix $\mathcal{I}_n (F_1, F_2, \ldots, F_n |\theta)$ associated with $F_1, F_2, \ldots, F_n.$
\end{theorem}
\vskip0.2in
\begin{proof}
Immediate.
\end{proof}

\vskip0.2in
 \begin{corollary}\label{cor3.1}
Let $X_i,$ $ 1 \leqslant i \leqslant m,$ $F_j,$ $1 \leqslant j \leqslant n$ be as in Theorem \ref{thm3.1}. Then
\begin{eqnarray*}
\lefteqn{\Lambda_{mn}(\theta) \mathcal{I}_n^{-} (F_1, F_2, \ldots, F_n | \theta) \Lambda_{mn} (\theta)^{\prime}}\\
&\geqslant& \Lambda_{mn-1} (\theta) \mathcal{I}_{n-1}^{-} (F_1, F_2, \ldots, F_{n-1} | \theta) \Lambda_{mn-1} (\theta)^{\prime}, \theta \in \Gamma  
\end{eqnarray*}
for $n \geqslant 2.$
 \end{corollary}
\vskip0.2in
\begin{proof}
This is immediate from the fact that both the sides of the inequality above are arrived at by taking supremum over certain sets in $\mathbb{R}^n,$ the set for the left hand side being larger than the set for the right hand side.
\end{proof}

We call the right hand side of the inequality in Theorem \ref{thm3.1} the CRB lower bound.
\vskip0.2in
\noindent{\bf Remark 1 \quad} Theorem \ref{thm3.1} and Corollary \ref{cor3.1}  imply the possibility of improving the CRB lower bound by searching for a larger class of $\mathcal{A}(\Gamma)$-linearly independent Fisher maps for a parametric family of states.

\vskip0.2in
\noindent{\bf Remark 2 \quad} The CRB lower bound in Theorem \ref{thm3.1} 
has some natural invariance properties. If $f_i,$ $1 \leqslant i \leqslant m$ are fixed and $X_i,$ $F_i (\theta), \rho(\theta)$ are changed respectively to $UX_iU^{\dagger},$ $U F_i(\theta) U^{\dagger},$ $U \rho (\theta)U^{\dagger}$ by a fixed unitary operator $U$ in $\mathcal{H}$ then the CRB lower bound in Theorem \ref{thm3.1} remains the same.

If the Fisher maps $F_j$ are replaced by
\begin{equation}
G_j (\theta) = \sum_{r=1}^n \alpha_{jr}(\theta) F_r (\theta), \quad 1 \leqslant j \leqslant n \label{eq3.8}
\end{equation}
where the matrix $A(\theta) = ((\alpha_{rs} (\theta)))$ is invertible for all $\theta$ then
$$\Lambda_{mn}(\theta) \mathcal{I}_n^{-} (F_1, \ldots, F_n | \theta) \Lambda_{mn}(\theta)^{\prime} =  \tilde{\Lambda}_{mn} (\theta) \mathcal{I}_n^{-} (G_1, \ldots, G_n|\theta) \tilde{\Lambda}_{mn} (\theta)^{\prime},$$
the tilde over $\Lambda_{mn}$ on the right hand side indicating that $G_i$'s are used in place of $F_i$'s.   In other words the CRB bound is invariant under $\mathcal{A}(\Gamma)$-linear invertible transformations of the form \eqref{eq3.8}.
\vskip0.2in
\noindent{\bf Example 3.1}
Let $\Gamma = (a,b),$ $\mathcal{H} = \mathbb{C}^n$ and let
$$\rho (\theta) = \diag \left (p_1(\theta), p_2(\theta), \ldots, p_n(\theta) \right ), \quad \theta \in \Gamma$$
be states in $\mathbb{C}_n$ with respect to the standard orthonormal basis, $\diag$  denoting diagonal matrix. An estimable function $f$ on $\Gamma$ has the form
$$f(\theta) = \sum_{i=1}^n \, a_i \, p_i (\theta)$$
where $a_i$ are real scalars. An unbiased estimator $X$ for $f$ is
$$X = \diag (a_1, a_2, \ldots, a_n).$$
Note that $p_i(\theta) \geqslant 0$ and $\sum\limits_i p_i (\theta) = 1$ $\forall$ $\theta \in \Gamma.$ Assume that $p_i (\theta)$ are differentiable in $\theta$ and $p_i (\theta) > 0$ $\forall$ $i, \theta.$ Then
$$F(\theta) = \diag \left (\frac{p_1^{\prime}(\theta)}{p_1(\theta)}, \frac{p_2^{\prime}(\theta)}{p_2(\theta)}, \ldots,   \frac{p_n^{\prime}(\theta)}{p_n(\theta)}  \right )$$
yields a Fisher map with
$$\mathcal{I} (F, F)(\theta) = \sum_{i=1}^n \,\,\frac{p_i^{\prime}(\theta)^2}{p_i(\theta)} $$
and 
$$\lambda (f, F) = \sum_{i=1}^n \,a_i \, p_i^{\prime}(\theta) = f^{\prime}(\theta). $$
Theorem \ref{thm3.1} for the single observable $X$ and single Fisher map yields
$$\Var \, \left ( Y | \theta \right ) \geqslant \,\frac{\left (\sum\limits_{i=1}^n \,a_i \, p_i^{\prime} (\theta) \right )^2}{\sum\limits_{i=1}^n \,\frac{p_i^{\prime} (\theta)^2}{p_i (\theta)}} \quad \forall \,\, \theta \in (a,b) $$
and any unbiased estimator $Y$ of $f.$ This is the Cram\'er-Rao inequality for finite sample spaces in classical mathematical statistics.

\vskip0.2in
\noindent{\bf Example 3.2}(Quantum version of Barankin's example \cite{ewb}, \cite{hlvt}).
Let $\rho(\theta)$ be an invertible density operator for every $\theta$ in $\Gamma.$ For any $\gamma \in \Gamma$ define
$$F_{\gamma} (\theta) = \rho (\gamma) \rho (\theta)^{-1} - 1. $$
Then $F_{\gamma}$ is a Fisher map and for any estimable function $f \in \mathcal{E}(\Gamma)$ we have
$$\lambda (f, F_{\gamma}) (\theta) = f(\gamma) - f(\theta). $$
The Fisher information form $\mathcal{I}$ satisfies
$$\mathcal{I} (F_{\gamma_{1}},  F_{\gamma_{2}} )(\theta) = \,\Re \,\tr \,\rho(\gamma_1) \rho(\theta)^{-1} \rho (\gamma_2) - 1 $$
If $X$ is an unbiased estimate of $f \in \mathcal{E}(\Gamma)$ one obtains as a special case the CRB bound
\begin{gather*}
\Var (X|\theta) \geqslant \left (f(\gamma_1)-f(\theta), f(\gamma_2)-f(\theta), \ldots, f(\gamma_n)-f(\theta) \right )\\
 \mathcal{I}^{-}_n (\gamma_1, \gamma_2, \ldots, \gamma_n, \theta) \left (f(\gamma_1)-f(\theta), \ldots, f(\gamma_n)-f(\theta) \right )^{\prime}
\end{gather*}
where $\mathcal{I}^{-}_n (\gamma_1, \gamma_2, \ldots, \gamma_n \theta)$ is the generalized inverse of the information matrix
$$\left ( \left (\Re \,\tr \, \rho (\gamma_i) \rho(\theta)^{-1} \rho (\gamma_j) - 1 \right ) \right ) $$
for any $\gamma_1, \gamma_2, \ldots, \gamma_n \in \Gamma.$

\vskip0.2in
\noindent{\bf Example 3.3}(Quantum Bhattacharya bound \cite{ab}).
Let $\Gamma \subseteq \mathbb{R}^d$ be a connected open set and let $\rho(\btheta),$ $\btheta \in \Gamma$ be a family of invertible states such that the correspondence $\btheta \rightarrow \rho(\btheta)$ is $C^m$-smooth. then every estimable function $f$ is also $C^m$-smooth. For any linear differential operator $D$ on $\Gamma$ with $C^m$-coefficients satisfying $D\,1 = 0$ define
$$F_D (\btheta) = (D \rho) (\btheta) \rho (\btheta)^{-1} $$
where $D$ is applied to every matrix entry of $\rho(\cdot)$ on the right hand side in some fixed orthonormal basis. Then $F_D$ is a Fisher map and the CRB tensor $\lambda$ satisfies
$$\lambda (f, F_D) (\btheta) = (Df) (\btheta) \quad \forall \,\, f \in \mathcal{E} (\Gamma).$$
If $D_1, D_2$ are two linear differential operators in $\Gamma$ with $C^m$-coefficients annihilating the constant function $1$ the Fisher information satisfies
$$\mathcal{I} (F_{D_{1}}, F_{D_{2}}) (\btheta) = \,\Re \,\,\tr \,\, (D_1 \rho) (\btheta) \,\rho(\btheta)^{-1} (D_2 \rho) (\btheta), \quad \btheta \in \Gamma.$$
If $X$ is an unbiased estimate of $f$ and $D_i,$ $1 \leqslant i \leqslant n$ are $C^m$-differential operators on $\Gamma$ then the CRB inequality has the form
$$\Var (X|\theta) \geqslant (D_1 f, \ldots, D_n f) (\btheta) \mathcal{I}_n^{-} (D_1, D_2, \ldots, D_n | \btheta) (D_1 f, \ldots, D_n f)(\btheta)^{\prime} $$
where $\mathcal{I}_n^{-} (D_1, D_2, \ldots, D_n | \btheta)$ is the generalized inverse of the positive semidefinite matrix
$$\left ( \left ( \Re \,\,\tr\,\, (D_i \rho) (\btheta) \rho(\btheta)^{-1} (D_j \rho)(\btheta) \right ) \right ), \quad i, j \in \{ 1,2, \ldots, n\}.$$

\vskip0.2in
\noindent{\bf Example 3.4 \,\,} 
Example 3.2 leads us to the following natural abstraction. Suppose $\Gamma$ is a $d$-dimensional $C^m$-manifold and $\btheta \rightarrow \rho (\btheta)$ is a $C^m$-smooth parametrization of states in $\mathcal{H}$ as $\btheta$ varies in $\Gamma.$ If $L$ is a smooth vector field on $\Gamma$ then
$$F_L (\theta) = (L \rho) (\btheta) \rho(\btheta)^{-1}, \quad \theta \in \Gamma $$
is a $C^m$-smooth Fisher map with respect to $\{\rho (\theta), \theta \in \Gamma \}$ under the assumption that $\rho (\theta)^{-1}$ exists for every $\theta.$ $C^m$-smooth Fisher maps constitute a $C^m(\Gamma)$-module and $\mathcal{E}(\Gamma) \subset C^m (\Gamma).$ The CRB tensor $\lambda$ and the Fisher information form $\mathcal{I}$ satisfy the relations
\begin{eqnarray*}
\lambda (f, F_L) (\theta) &=& (Lf) (\theta) \\
\mathcal{I}(F_L, F_M) (\theta) &=& \Re \,\tr \, (L \rho) (\theta) \rho(\theta)^{-1} (M\rho) (\theta)
\end{eqnarray*}
for any two vector fields $L,M.$ As a special case of the CRB inequality we have for any unbiased estimator $X$ of $f \in \mathcal{E}(\Gamma),$
$$\Var (X| \theta) \geqslant \frac{(Lf)(\theta)^2}{\tr \,\,\rho(\theta)^{-1} (L\rho)(\theta)^2}, \quad \theta \in \Gamma$$
for any $C^m$-smooth vector field $L.$

As a special case of the example above, consider a connected Lie group $\Gamma$ with Lie algebra $\mathcal{G}.$ Let 
$$\rho (g) = U_g \,\,\rho_0 \,\,U_g^{\dagger}, \quad g \in \Gamma$$
where $\rho_0$ is a fixed invertible state. Any element $L$ of $\mathcal{G}$ is looked upon as a left invariant vector field on $\Gamma.$ let $U_{\exp \,\,t\,\,L} = \exp \,\,t\,\,\pi(L),$ $t \in \mathbb{R},$ $L \in \mathcal{G}$ where $L \rightarrow \pi (L)$ is a representation of $\mathcal{G}$ in $\mathcal{H}.$ Then the CRB inequality yields
\begin{equation}
\Var (X|g) \geqslant \frac{((Lf)(g))^2}{\tr \,\,\rho_0^{-1} \left [ \pi (L), \rho_0 \right ]^2}  \quad \forall \quad L \in \mathcal{G} \label{eq3.9}
\end{equation}
where $X$ is an unbiased estimator of $f.$ If $L_i,$ $1 \leqslant i \leqslant d$ is a basis for $\mathcal{G}$ and the nonnegative definite matrix $\mathcal{I}_d$ is defined by
$$\mathcal{I}_d = \left ( \left ( \Re \,\,\tr\,\,\rho_0^{-1} \left [\pi (L_i), \rho_0 \right ] \left [\pi (L_j), \rho_0 \right ] \right ) \right ), \quad i, j \in \{1,2,\ldots,d\} $$
then a maximization over all elements $L$ in $\mathcal{G}$ on the right hand side of \eqref{eq3.9} yields
$$\Var (X|g) \geqslant (L_1 \,f, L_2 \, f, \ldots, L_d f)(g) \mathcal{I}_d^{-} (L_1 f, L_2 f, \ldots, L_d f) (g)^{\prime},$$
$\mathcal{I}_d^{-}$ being the generalized inverse of $\mathcal{I}_d.$

\vskip0.2in
\noindent{\bf Example 3.5}(adapted from \cite{slb}).
We now consider an example in which the different states $\rho(\theta)$ may fail to have an inverse, indeed, their ranges need not be the same. Let $\Gamma \subseteq \mathbb{R}^d$ be an open domain and let $\rho(\btheta),$ $\btheta \in \Gamma$ obey the set of linear partial differential equations of the form
\begin{equation}
\frac{\partial \rho}{\partial \theta_j} = \frac{1}{2} \left (L_j (\btheta) \,\rho (\btheta) + \rho (\btheta) L_j (\btheta)^{\dagger} \right ), \quad 1 \leqslant j \leqslant d \label{eq3.10} 
\end{equation}
where the operators $L_j (\btheta) \in \mathcal{B}(\mathcal{H}).$ Taking trace on both sides we see that
$$\Re \,\,\tr\,\, \rho(\btheta) L_j (\btheta) = 0, \quad 1 \leqslant j \leqslant d, \quad \btheta \in \Gamma.$$ 
If $\mathcal{I}m \,\,\tr\,\, \rho(\btheta) L_j(\btheta) = m_j (\btheta)$ we can replace in \eqref{eq3.10} $L_j (\btheta)$ by $L_j (\btheta) - im_j(\btheta)$ without altering the differential equations. Hence we may assume, without loss of generality, that in \eqref{eq3.10}
\begin{equation}
\tr \,\, \rho(\btheta) L_j (\btheta) = 0, \quad 1 \leqslant j \leqslant d, \quad \btheta \in \Gamma.\label{eq3.11}
\end{equation}
We then say that the states $\rho(\btheta)$ which obey \eqref{eq3.10} and \eqref{eq3.11} constitute a {\it Liapunov family.}

A special case of such a Liapunov family of states is obtained when $d=1$ and
$$\rho(\theta) = p (\theta) e^{\frac{1}{2} \theta L} \, \rho_0 \,e^{\frac{1}{2} \theta L^{\dagger}}, \quad \theta \in \mathbb{R} $$
where $L$ is a fixed operator in $\mathcal{H},$ $\rho_0$ is a fixed state and
$$p(\theta) = \left \{\tr \,\, \rho_0 \,\,e^{\frac{1}{2} \theta L^{\dagger}} \,\,e^{\frac{1}{2} \theta L} \right \}^{-1}. $$
Then
$$\rho^{\prime}(\theta) = \frac{1}{2} \left \{ \left ( \frac{p^{\prime}(\theta)}{p(\theta)} + L\right ) \rho(\theta) + \rho(\theta) \left (\frac{p^{\prime}(\theta)}{p(\theta)} + L \right )^{\dagger} \right \}. $$
If $\rho_0 = | \psi_0 \rangle \langle \psi_0 |$ is a pure state then every $\rho(\theta)$ is a pure state. Thus ${\rm rank}\,\rho(\theta) = 1 \,\,\forall \,\,\theta \in \mathbb{R}$ and we have a situation where $\{ \rho (\theta) \}$ admits a `score operator function' with a classical part $p^{\prime}/p$ and a quantum part $L.$

Going back to the Liapunov family satisfying \eqref{eq3.10} and \eqref{eq3.11} we observe that each of the maps $\btheta \rightarrow L_j (\btheta),$ $1 \leqslant j \leqslant d$ is a Fisher map. Indeed, if $X$ is a balanced observable we have
\begin{eqnarray*}
 0 &=& \frac{\partial}{\partial \theta_j} (\tr \, \rho(\btheta) X) \\
&=& \frac{1}{2} \,\tr \, \left (L_j (\btheta) \rho(\btheta) X + \rho (\btheta) L_j (\btheta)^{\dagger} X \right ) \\
&=& \frac{1}{2}\,\tr \,\rho (\theta) \left ( L_j (\theta)^{\dagger} X + X L_j (\btheta) \right ).
\end{eqnarray*}
For any estimable function $f$
$$\lambda (f, L_j) (\btheta) = \frac{\partial f}{\partial \theta_j}$$
and the Fisher information form $\mathcal{I}$ satisfies
$$\mathcal{I} (L_i, L_j) (\btheta) = \,\Re \,\, \tr \,\, \rho (\btheta) L_i (\btheta)^{\dagger} L_j (\btheta).$$
If we write
$$\mathcal{I}_d (\btheta) = \left ( \left (\mathcal{I} (L_i, L_j)(\btheta) \right ) \right ), \quad i, j \in \{ 1,2,\ldots, d\} $$
then the CRB inequality assumes the form
$$\Var (X|\btheta) \geqslant (\nabla f) (\btheta) \mathcal{I}_d^{-} (\btheta) (\nabla f)(\btheta)^{\prime}$$
for any unbiased estimator $X$ of $f,$ $\nabla f$ being the gradient vector $\left (\frac{\partial f}{\partial \theta}, \frac{\partial f}{\partial \theta_2}, \ldots, \frac{\partial f}{\partial \theta_d} \right ).$

In the special case $d=1$ introduced in the course of the discussion above the CRB bound assumes the form
$$\Var (X|\theta) \geqslant \frac{\left (f^{\prime} (\theta) \right )^2}{\tr \,\, \rho(\theta) \left (\frac{p^{\prime} (\theta)}{p(\theta)} + L  \right )^{\dagger}      \left (\frac{p^{\prime} (\theta)}{p(\theta)} + L  \right )     }. $$

If $\rho(\btheta),$ $\sigma(\btheta),$ $\btheta \in \Gamma$ are Liapunov families of states in Hilbert spaces $\mathcal{H},$ $\mathcal{K}$ respectively with coefficients $L_j(\btheta),$ $M_j(\btheta)$ in the respective differential equations corresponding to \eqref{eq3.10} then the tensor product states $\rho(\btheta) \otimes \sigma(\btheta),$ $\btheta \in \Gamma$ constitute again a Liapunov family with the coefficients $L_j (\btheta) \otimes 1 + 1 \otimes M_j (\theta),$ $1 \leqslant j \leqslant d$ in the differential equations corresponding to \eqref{eq3.10} and its Fisher information form satisfies
\begin{eqnarray*}
\lefteqn{\mathcal{I} \left (L_i \otimes 1 + 1 \otimes M_i, L_j \otimes 1 + 1 \otimes M_j \right ) (\btheta)} \\
&=& \mathcal{I} (L_i, L_j) (\btheta) + \mathcal{I} (M_i, M_j) (\btheta).
\end{eqnarray*}
\vskip0.2in
\noindent{\bf Eexample 3.6 \,}
Our last example is the case when $\rho(\btheta)$ is a mixture of the form
$$\rho (\btheta)  = \sum_{r=1}^N \,p_r (\btheta) \,\rho_r (\btheta)$$
where $\{p_r (\btheta), 1 \leqslant r \leqslant N \}$ is a family of probability distributions on the finite set $\{1,2,\ldots, N\}$ indexed by $\btheta \in \Gamma$ and for each fixed $r,$ $\{\rho_r (\btheta), \btheta \in \Gamma \}$ is a Liapunov family of states obeying the differential equations
$$\frac{\partial \rho_r}{\partial \theta_j} = \frac{1}{2} \left \{L_{rj} (\btheta) \,\rho_r (\btheta) + \rho_r (\btheta) L_{rj}(\btheta)^{\dagger} \right \}, \quad 1 \leqslant j \leqslant d, \,\,\btheta \in \Gamma $$
and the conditions
$$\tr \,\rho_r (\btheta) \,L_{rj} (\btheta) = 0 \quad \forall \quad \btheta \in \Gamma.$$
Let now $f_i, 1 \leqslant i \leqslant m$ be estimable functions with respect to $\{\rho (\btheta), \btheta \in \Gamma\}$ and let $X_i$ be any unbiased estimator of $f_i$ for each $i.$ Differentiating with respect to $\theta_j$ the identity
$$\tr \, \rho(\btheta) (X_i - f_i (\theta)) = 0$$
we get
\begin{equation}
\frac{\partial f_i}{\partial \theta_j} = \sum_{r=1}^N \,p_r (\btheta) \,\Re \,\tr \,M_{rj} (\btheta) \rho_r (\btheta) (X_i - f_i (\btheta)) \label{eq3.12}
\end{equation}
where
\begin{equation}
M_{rj} (\btheta) = p_r (\btheta)^{-1} \frac{\partial p_r}{\partial \theta_j} + L_{rj} (\btheta).      \label{eq3.13}
\end{equation}
Multiplying both sides of  \eqref{eq3.12} by real scalars $a_i b_j$ and adding over $i$ and $j$ we get
$$\mathbf{a}^{\prime} \left ( \left ( \frac{\partial f_i}{\partial \theta_j} \right ) \right ) \mathbf{b} = \sum_{r=1}^N \,p_r (\btheta) \,\tr\, \left (\sum_{j=1}^d b_j M_{rj} (\btheta) \right ) \rho_r (\btheta) \left (\sum_{i=1}^m a_i (X_i - f_i (\btheta)) \right ). $$
Applying Cauchy-Schwarz inequality to each trace scalar product on the right hand side followed by the same inequality to the scalar product with respect to the probability distribution $p_1(\btheta), p_2(\btheta), \ldots, p_N(\btheta)$ we obtain
\begin{eqnarray}
\left (\mathbf{a}^{\prime} \left ( \left ( \frac{\partial f_i}{\partial \theta_j} \right ) \right ) \,\mathbf{b} \right )^2 \leqslant \left \{\sum_{r=1}^N p_r (\btheta) \,\tr\, \left ( \sum_{j=1}^d b_j M_{rj} (\btheta) \right )  \right . \nonumber\\
\left . \rho_r (\btheta) \left (\sum_{j=1}^d b_j M_{rj}(\btheta) \right )^{\dagger} \right \} \mathbf{a}^{\prime}  \,\Cov (X_1, \ldots, X_m | \btheta)\, \mathbf{a}      \label{eq3.14}
\end{eqnarray}
Let
\begin{eqnarray*}
\Psi_r (\btheta) &=& \left ( \left ( \Re \,\,\tr \, \, \rho_r (\btheta) M_{ri} (\btheta)^{\dagger} M_{rj} (\btheta) \right )\right ), \quad i, j \in \{ 1,2, \ldots, d \}, \\
\Psi(\btheta) &=& \sum_{r=1}^N p_r (\btheta) \Psi_r (\btheta).
\end{eqnarray*}
Then the validity of \eqref{eq3.14} for all $a_i, b_j,$ $1 \leqslant i \leqslant m,$ $1 \leqslant j \leqslant d$ implies
$$\Cov (X_1, X_2, \ldots, X_m|\btheta) \geqslant \left ( \left ( \frac{\partial f_i}{\partial \theta_j} \right ) \right ) \Psi^{-} (\btheta) \left ( \left (\frac{\partial f_i}{\partial \theta_j} \right ) \right )^{\prime},$$
the super index - in $\Psi$ indicating its generalized inverse.

\section{Estimators based on generalized measurements}
\setcounter{equation}{0}

As in Section 3 we consider a parametric family $\{\rho (\theta), \theta \in \Gamma\}$ of states of a finite level quantum system in a Hilbert space $\mathcal{H}$ and a real-valued parametric function $f$ on $\Gamma.$ In order to estimate $f$ we now look at a generalized measurement $\mathcal{L} = (S,L)$ as described in Definition \ref{def2.1}. Choose a real-valued function $\varphi$ on $S$ and if the outcome of $\mathcal{L}$ is $s$ then evaluate $\varphi(s)$ and treat it as an estimate of $f(\theta).$ We say that $(\mathcal{L}, \varphi)$ is an {\it unbiased estimator} of $f$ if

\begin{equation}
\sum_{s \in S} \varphi(s) \,\tr\, \rho(\theta) L(s)^{\dagger} L(s) = f(\theta) \quad \forall \,\,\theta \in \Gamma. \label{eq4.1}
\end{equation}
Indeed, it may be recalled from Section 2 that $\tr \,\,\rho(\theta) L(s)^{\dagger} L(s)$ is the probability of the outcome $s$ if the unknown parameter is $\theta.$ Then the variance of $(\mathcal{L}, \varphi)$ is given by
\begin{equation}
\Var (\mathcal{L}, \varphi | \theta) = \sum_{s \in S} \,\varphi(s)^2 \,\,\tr\,\, \rho(\theta) L(s)^{\dagger} L(s) - f (\theta)^2. \label{eq4.2}
\end{equation}
If we write
\begin{equation}
X = \sum_{s \in S} \,\varphi (s) L(s)^{\dagger} L(s)  \label{eq4.3}
\end{equation}
Then $X$ is an observable and \eqref{eq4.1} shows that $X$ is an unbiased estimator of $f$ whenever $(\mathcal{L}, \varphi)$ is an unbiased estimator of $f.$ However, $\Var(X|\theta)$ need not be the same as $\Var (\mathcal{L}, \varphi| \theta).$

In \eqref{eq4.1} put $T(s) = L(s)^{\dagger} L(s), s \in S.$ Then $T(s) \geqslant 0$ and by Definition \ref{def2.1}, $\sum\limits_{s \in S} T(s) = I.$ In other words $\{T(s), s \in S \}$ is a positive operator-valued distribution on $S$ with total operator mass $I.$ By a well-known theorem of Naimark \cite{ash}, \cite{KRP} one can imbed the Hilbert space $\mathcal{H}$ isometrically in a larger Hilbert space $\widehat{\mathcal{H}} = \mathcal{H} \otimes \mathcal{K}$ with $\dim \,\mathcal{K} < \infty$ and construct mutually orthogonal projection operators on $\widehat{\mathcal{H}}$ with the block operator form
\begin{equation}
E(s) = \left [\begin{array}{c|c} T(s) & M(s) \\ \hline M(s)^{\dagger} & N(s)  \end{array} \right ], \quad s \in S \label{eq4.4}
\end{equation}
satisfying the following:
\begin{itemize}
 \item[(i)] $\sum\limits_{s \in S} \,\,E(s) = I,$
\item[(ii)] $\left \{ E(s) \left [\begin{array}{c} u \\ 0 \end{array} \right ], \,\,s \in S, \,\, u \in \mathcal{H} \right \} \quad {\rm spans} \,\,\widehat{\mathcal{H}}.$
\end{itemize}
Such a dilation of $T(\cdot)$ in $\mathcal{H}$ to $E(\cdot)$ in $\widehat{\mathcal{H}}$ is unique upto a natural unitary isomorphism.

Now we go back to the unbiased estimator $(\mathcal{L}, \varphi)$ of $f$ described in \eqref{eq4.1}. Put
\begin{eqnarray*}
\widehat{\rho} (\theta) &=& \left [ \begin{array}{c|c} \rho(\theta) & 0 \\ \hline 0 & 0 \end{array} \right ], \\
\widehat{X} &=& \sum_{s \in S} \,\,\varphi (s) \,\,E(s).
\end{eqnarray*}
Then $\{\widehat{\rho} (\theta), \theta \in \Gamma\}$ is a parametric family of states in $\widehat{\mathcal{H}},$ $\widehat{X}$ is an observable in $\widehat{\mathcal{H}}$ and equations \eqref{eq4.1} and \eqref{eq4.4} imply that $\tr \,\, \widehat{\rho}(\theta) \widehat{X} = f (\theta).$ Furthermore
\begin{eqnarray*}
\Var (\widehat{X}|\theta) &=& \tr \,\,\widehat{\rho} (\theta) (\widehat{X} - f(\theta))^2 \\
&=& \sum_{s \in S} \,\varphi(s)^2 \,\tr \,\, \rho (\theta) T(s) - f (\theta)^2 \\
&=& \Var(\mathcal{L}, \varphi | \theta).
\end{eqnarray*}
Thus $\widehat{X}$ is an unbiased estimator of $f$ with respect to $\{\widehat{\rho} (\theta), \theta \in \Gamma \}$ with the same variance as the unbiased estimator $(\mathcal{L}, \varphi)$ based on generalized measurement for the original family of states.

If $F$ is a Fisher map for $\{\rho(\theta), \theta \in \Gamma \}$ then $\widehat{F}$ defined by
$$\widehat{F} (\theta) = \left [ \begin{array}{c|c} F(\theta) & 0 \\ \hline 0 & 0 \end{array} \right ], \quad \theta \in \Gamma $$
is a Fisher map for $\{\widehat{\rho} (\theta), \theta \in \Gamma \}$ in $\widehat{\mathcal{H}}.$ If $\widehat{\mathcal{I}}$ is the Fisher information form for $\{\widehat{\rho} (\theta), \theta \in \Gamma \}$ we have
$$ \widehat{\mathcal{I}} (\widehat{F}_1, \widehat{F}_2)(\theta) = \mathcal{I} (F_1, F_2) (\theta).$$
Thus from Theorem \ref{thm3.1} we conclude the following theorem.
\vskip0.1in
\begin{theorem}\label{thm4.1}
Let $\{\rho (\theta), \theta \in \Gamma \}$ be a parametric family of states of a finite level quantum system in a Hilbert space $\mathcal{H}$ and let $(\mathcal{L}, \varphi)$ be any unbiased estimator of a real-valued parametric function $f$ based on a generalized measurement $\mathcal{L}$ and a real scalar function $\varphi$ on the set of values of the measurement. Suppose $F_j, 1 \leqslant j \leqslant n$ are Fisher maps for $\{\rho (\theta), \theta \in \Gamma \}.$ Then
\begin{eqnarray*}
\lefteqn{\Var \left ( (\mathcal{L}, \varphi)| \theta \right ) \geqslant \left (\lambda(f, F_1), \lambda (f,F_2), \ldots, \lambda(f, F_n) \right )} \\ 
&& \mathcal{I}_n^{-} \left (\lambda(f,F_1), \lambda(f,F_2), \ldots, \lambda(f, F_n) \right )^{\prime} (\theta)
\end{eqnarray*}
where $\lambda$ is the CRB tensor and $\mathcal{I}_n^{-}$ is the generalized inverse of the information matrix
$$\mathcal{I}_n = \left ( \left ( \mathcal{I} (F_i, F_j) \,\right ) \right ), \quad i, j \in \{ 1,2, \ldots, n \}.$$
\end{theorem}
\vskip0.1in
\begin{proof}
Immediate.
\end{proof}

We shall briefly consider the case of estimating many parametric functions $f_i (\theta),$ $1 \leqslant i \leqslant m.$ In order to estimate them it appears that several generalized measurements are to be made. Such measurements have to be made in succession. As directed in Section 2 we may treat them all as a single compound generalized measurement $\mathcal{L} = (L, S).$ Let $(\mathcal{L}, \varphi_i)$ be an unbiased estimator of $f_i$ for each $i.$ Thus the measurement $\mathcal{L}$ is carried out and if the outcome is $s \in S$ then $\varphi_i (s)$ is the estimate of $f_i(\theta).$ The probability for the outcome $s$ is $\tr \,\, \rho(\theta) L(s)^{\dagger} L(s).$ Thus the covariance matrix of the different estimators is given by
\begin{eqnarray}
\lefteqn{\Cov \left ( \mathcal{L}, \varphi_1, \varphi_2, \ldots, \varphi_m | \theta  \right )           } \nonumber \\
&=&  \left ( \left ( \tr \,\, \rho(\theta) \sum_{s \in S} \varphi_i(s) \varphi_j(s) L(s)^{\dagger}  L(s) - f_i(\theta) f_j (\theta) \right ) \right ), \nonumber \\
&&  \quad i. j \in \{1,2, \ldots, m\}. \label{eq4.5}
\end{eqnarray}
As in the discussion preceding Theorem \ref{thm4.1} we can construct the Naimark dilation $\{ E(s), s \in S \}$ for the positive operator-valued distribution $\{L(s)^{\dagger} L(s), s \in S \}$ in an enlarged Hilbert space and view the covariance matrix \eqref{eq4.5} as
$$\Cov \left (\widehat{X}_1,  \widehat{X}_2, \ldots, \widehat{X}_m | \theta \right )$$
for the observables $\widehat{X}_i, = \sum\limits_s \,\varphi_i(s) E(s)$ with respect to the states $\widehat{\rho}(\theta).$ This at once leads us to the CRB matrix inequality
\begin{eqnarray*}
\lefteqn{\Cov \left ( \mathcal{L}, \varphi_1, \varphi_2, \ldots, \varphi_m|\theta \right ) \geqslant \left ( \left (\lambda(f_i, F_j) \right ) \right ) \left ( \left (\mathcal{I}_n^{-} (F_p, F_q) \right ) \right ) }\\
&& \left ( \left ( \lambda (f_i, F_j) \right ) \right )^{\prime} (\theta), \,\, 1\leqslant i \leqslant m; \,\,j,p,q \in \{ 1,2, \ldots, n\}. 
\end{eqnarray*}
for any set $\{F_j, 1 \leqslant j \leqslant n \}$ of Fisher maps, $\lambda$ being the CRB tensor, $\left ( \left ( \mathcal{I}_n (F_i, F_j) \right ) \right )$ the Fisher information matrix with respect to $\{F_j, 1 \leqslant j \leqslant n \}$ and the super index - denoting generalized inverse.

\vskip0.2in
\noindent{\it Acknowledgement \quad} The author thanks H. Parthasarathy for several useful comments and also pointing out the references \cite{jm} and \cite{hlvt} in the physics and engineering literature. He thanks B. V. Rao for bringing his attention to the very rich survey article \cite{dct}.

\end{document}